# Optimal Joint Bidding and Pricing of Profit-seeking Load Serving Entity

Hanchen Xu, *Student Member, IEEE,* Kaiqing Zhang, *Student Member, IEEE,* and Junbo Zhang, *Member, IEEE*

***Abstract*—The demand response provides an opportunity for load serving entities (LSEs) that operate retail electricity markets (REMs) to strategically purchase energy and provide reserves in wholesale electricity markets (WEMs). This paper concerns with the problem of simultaneously determining the optimal energy bids and reserve offers an LSE submits to the WEM as well as the optimal energy and reserve prices it sets in the REM so as to maximize its profit. To this end, we explicitly model the tri-layer market structure that consists of a WEM, a REM, and a set of end user customers, so as to capture the coupling between the bidding problem and the pricing problem. Based on the tri-layer market model, we then formulate the joint bidding and pricing problem as a bi-level programming problem and further transform it into a single-level mixed integer linear programming problem, which can be solved efficiently. Numerical studies using the IEEE test cases are presented to illustrate the application of the proposed methodology as well as to reveal several interesting characteristics of the LSE's profit-seeking behavior.**

## NOMENCLATURE

*Acronyms*:

| | | |
|---|---|---|
| DR | : | demand response |
| EUC | : | end user customer |
| GEN | : | generator |
| ISF | : | injection shifting factor |
| ISO | : | independent system operator |
| LF | : | loss factor |
| LMP | : | locational marginal price |
| LSE | : | load serving entity |
| MILP | : | mixed integer linear programing |
| REM | : | retail electricity market |
| WEM | : | wholesale electricity market |

*Symbols*:

| | | |
|---|---|---|
| $\mathcal{N}$ | : | index set of buses |
| $\mathcal{L}$ | : | index set of lines |
| $\mathcal{G}$ | : | index set of GEN offers |
| $\mathcal{D}$ | : | index set of LSE bids |
| $\tilde{\mathcal{N}}$ | : | index set of buses where the LSE submits bids |
| $\tilde{\mathcal{D}}$ | : | index set of the LSE's bids |
| $\tilde{\mathcal{D}}_k$ | : | index set of the LSE's bids submitted at bus $i_k \in \tilde{\mathcal{N}}$ |
| $\mathcal{C}$ | : | index set of EUCs associated with the LSE |
| $\mathcal{C}_k$ | : | index set of EUCs associated with bus $i_k \in \tilde{\mathcal{N}}$ |
| $M$ | : | a sufficiently large number |
| $\overline{\boldsymbol{f}}$ | : | maximum line active power flow vector (MW) |
| $\boldsymbol{\Lambda}$ | : | LF vector |
| $\boldsymbol{\Psi}$ | : | ISF matrix |
| $\boldsymbol{a}^g$ | : | GEN marginal energy offer price vector (\$/MWh) |
| $\boldsymbol{b}^g$ | : | GEN marginal reserve offer price vector (\$/MW) |
| $\overline{\boldsymbol{p}}^g$ | : | GEN maximum active power supply vector (MW) |
| $\underline{\boldsymbol{p}}^g$ | : | GEN minimum active power supply vector (MW) |
| $\overline{\boldsymbol{r}}^g$ | : | GEN maximum reserve vector (MW) |
| $\boldsymbol{a}^d$ | : | LSE marginal energy bid price vector (\$/MWh) |
| $\boldsymbol{b}^d$ | : | LSE marginal reserve offer price vector (\$/MW) |
| $\overline{\boldsymbol{p}}^d$ | : | LSE maximum active power demand vector (MW) |
| $\underline{\boldsymbol{p}}^d$ | : | LSE minimum active power demand vector (MW) |
| $\overline{\boldsymbol{r}}^d$ | : | LSE maximum reserve vector (MW) |
| $\boldsymbol{\Phi}^g$ | : | GEN offer to bus mapping matrix |
| $\boldsymbol{\Phi}^d$ | : | LSE offer to bus mapping matrix |
| $\overline{\alpha}$ | : | maximum allowed energy price for EUCs (\$/MWh) |
| $\underline{\alpha}$ | : | minimum allowed energy price for EUCs (\$/MWh) |
| $\overline{\beta}$ | : | maximum allowed reserve price for EUCs (\$/MWh) |
| $\underline{\beta}$ | : | minimum allowed reserve price for EUCs (\$/MWh) |
| $\alpha_k$ | : | energy price for EUCs in $\mathcal{C}_k$ (\$/MWh) |
| $\beta_k$ | : | reserve price for EUCs in $\mathcal{C}_k$ (\$/MW) |
| $\boldsymbol{c}^{(k)}$ | : | energy benefit function coefficient vector of EUCs in $\mathcal{C}_k$ (\$/MWh) |
| $\boldsymbol{b}^{(k)}$ | : | reserve cost function coefficient vector of EUCs in $\mathcal{C}_k$ (\$/MW) |
| $\boldsymbol{p}$ | : | net active power injection vector (MW) |
| $\boldsymbol{p}^g$ | : | active power supply vector (MW) |
| $\boldsymbol{p}^d$ | : | active power demand vector (MW) |
| $\boldsymbol{r}^g$ | : | GEN reserve vector (MW) |
| $\boldsymbol{r}^d$ | : | LSE reserve vector (MW) |
| $\pi_i$ | : | LMP at bus $i \in \mathcal{N}$ (\$/MWh) |
| $\boldsymbol{x}^{(k)}$ | : | energy demand vector of EUCs in $\mathcal{C}_k$ (MW) |
| $\boldsymbol{y}^{(k)}$ | : | reserve provision vector of EUCs in $\mathcal{C}_k$ (MW) |
| $\mathbf{1}$ | : | all-ones vector with appropriate dimension |
| $\mathbf{0}$ | : | all-zeros vector with appropriate dimension |

H. Xu and K. Zhang are with the Department of Electrical and Computer Engineering, University of Illinois at Urbana-Champaign, Urbana, IL 61801, USA. E-mail: {hxu45, kzhang66}@illinois.edu.

J. Zhang is with the School of Electric Power, South China University of Technology, Guangzhou 510641, China. E-mail: epjbzhang@scut.edu.cn.

## I. INTRODUCTION

THE demand response (DR) has proven to be a key component of the smart grid paradigm [1], [2] and has played a critical role in maintaining power system reliability [3]. A paramount driving force behind the impressive success of the DR is the market mechanisms in both wholesale electricity markets (WEMs) operated by independent system operators (ISOs) and retail electricity markets (REMs) operated by load serving entities (LSEs) [3]. Conventionally, a LSE purchases energy in a WEM to meet the typically inelastic demands and charge a fixed tarrif from the end user customers (EUCs) that participate in its REM. The rapid development of various DR programs such as real-time pricing has significantly increased the demand-side flexibility, enabling more elasticity in the modeling of the load [4], [5]. With such flexibility, the LSEs can now embrace the opportunities to strategically purchase energy, and additionally, provide reserves in the WEM, possibly with the objective of maximizing its own profits [6], [7].

In most conventional market models such as the one in [8], the LSE acts as an aggregator for the EUCs and submits the energy bids and reserve offers to the ISO. However, in reality, the LSE as a middleman has the capability to determine a retail price different from the one posted by the ISO for its own profit [9]. Moreover, the LSE may strategically change its energy bids and reserve offers submitted to the WEM in order to impact market clearing results in some circumstances to increase its profit. As such, in the aforementioned market environment, a profit-seeking LSE is faced with two problems: (i) determining optimal energy bids and reserve offers it submits to the ISO, referred to as the bidding problem; and (ii) determining optimal energy and reserve prices it charges from or pays to the EUCs, referred to as the pricing problem.

Existing works are mostly only concerned with either one of the two problems. For example, despite of different models, [10], [11] both focus on the problem of optimal bidding, where [10] exploits the formulation of mathematical program with equilibrium constraints while [11] resorts to the modeling of supply function equilibrium. Moreover, [12], [13] consider developing the bidding strategy for specific LSEs, i.e., the aggregator of electric vehicles and unmanned aerial vehicles, both in day-ahead markets. Regrading pricing, [14] proposes a Stackelberg game between LSEs and EUCs to maximize the revenue of each LSE while [15] addresses the same problem under the assumption that the LSE is capable of learning/estimating the EUC's power consumption patterns for its own profits. Similarly, [16] proposes a dynamic pricing strategy also with the leader-follower game-theoretic model. In [17], an energy management system (EMS), which resembles the role of an LSE, coordinates the price-responsive demands to maximize their utility. In this regard, the EMS is not a profit-seeking entity as we consider in the present work.

Yet, from the perspective of an LSE, the optimal solutions to its bidding problem and pricing problem are indeed inherently coupled. Specifically, these two problems are coupled by the physical constraint that the LSE must balance the energy purchased in the WEM and that sold in the REM, as well as the reserve sold in the WEM and that purchased in the REM. Moreover, the total profit gained by the LSE depends on both the wholesale prices and the retail prices, as well as the cleared quantities of energy/reserve. In the meantime, the cleared quantities of energy/reserve and their prices in WEM/REM have direct impacts on each other. As a result, the separate consideration of these two problems will lead to a situation where the LSE fails to make the energy/reserve balance and cannot harvest the highest profit. Therefore, the joint consideration of the bidding problem and the pricing problem is indeed more desirable. However, to the best of our knowledge, the joint bidding and pricing problem has not been well studied yet. Indeed, the challenge of simultaneously determining the optimal bids and prices lies in formulating this problem in a solvable form.

To this end, we explicitly model the tri-layer market structure, which consists of a WEM, a REM, and a set of EUCs, so as to capture the structural characteristics of actual markets. This model, referred to as the tri-layer market model, allows us to explicitly take into account the inherent coupling between the bidding problem and the pricing problem. Based on the tri-layer market model, we then formulate the joint bidding and pricing problem as a bi-level programming problem, in which the profit maximization problem of the LSE serves as the upper-level problem, and the WEM clearing problem and the EUC benefit maximization problem serve as the lower-level problem. Yet, this bi-level programming problem, which has two optimization problems as constraints, is generally hard to solve [18]; therefore, we further prove that it can be transformed into a single-level mixed integer linear programming (MILP) problem, leveraging similar linearization techniques as in [19]. As a result, the joint bidding and pricing problem can be solved efficiently using existing commercial solvers. We show through numerical studies that the proposed model and formulation can simultaneously solve the bidding problem and the pricing problem. The impacts of such profit-seeking behavior on market efficiency are also explored numerically.

The major contributions of this paper are the following: (i) we formulate the joint bidding and pricing problem of an LSE as a bi-level programming problem based on the tri-layer market model; (ii) we transform the hard-to-solve bi-level programming problem, in which the WEM clearing problem and the EUC benefit maximization problem are embedded, to an equivalent single-level MILP problem and thus can be solved efficiently; and (iii) we validate the effectiveness of the proposed model via numerical simulations and reveal several interesting characteristics of the LSE's profit-seeking behavior. To the best of our knowledge, this is the first paper to study the joint bidding and pricing problem.

The remainder of this paper is organized as follows. Section II introduces the tri-layer market model, based on which the joint bidding and pricing problem is formulated as a bi-level programming problem in Section III. The bi-level programming problem is further transformed into a single-level MILP problem in Section IV. In Section V, the application of the proposed methodology is illustrated. Section VI concludes our work and discusses potential future work directions.

## II. Tri-layer Market Model

In this section, we introduce the tri-layer market model, the structure of which is illustrated in Fig. 1. The tri-layer market model consists of an upper layer representing a WEM, a middle layer representing a REM, and a lower layer representing a set of EUCs that participate in the REM.

### A. Upper Layer: Wholesale Electricity Market

Assume the ISO operates a transmission network that consists of a set of buses indexed by $\mathcal{N} = \{1, \cdots, N\}$ and a set of transmission lines indexed by $\mathcal{L} = \{1, \cdots, L\}$. Let $f_i$ denote the maximum active power flow on line $i$, $i \in \mathcal{L}$, and define $\overline{\boldsymbol{f}} = [f_1, \cdots, f_L]^\top$. Typically, in the WEM clearing model, the nonlinear relationship between the system losses and the active power injections can be approximated by the so-called loss factors (LFs), and that between the line flows and the active power injections can be approximated by the injection shifting factors (ISFs) [20], [21]. Let $\Lambda_i$ denote the LF at bus $i$, $i \in \mathcal{N}$, and define $\boldsymbol{\Lambda} = [\Lambda_1, \cdots, \Lambda_N]^\top$. Let $\boldsymbol{\Psi} \in \mathbb{R}^{L \times N}$ denote the ISF matrix, where $\Psi_{ij}$ is the sensitivity of the power flow on line $i$ with respect to net power injection at bus $j$. Denote the net power injections at bus $i$ by $p_i$ and define $\boldsymbol{p} = [p_1, \cdots, p_N]^\top$.

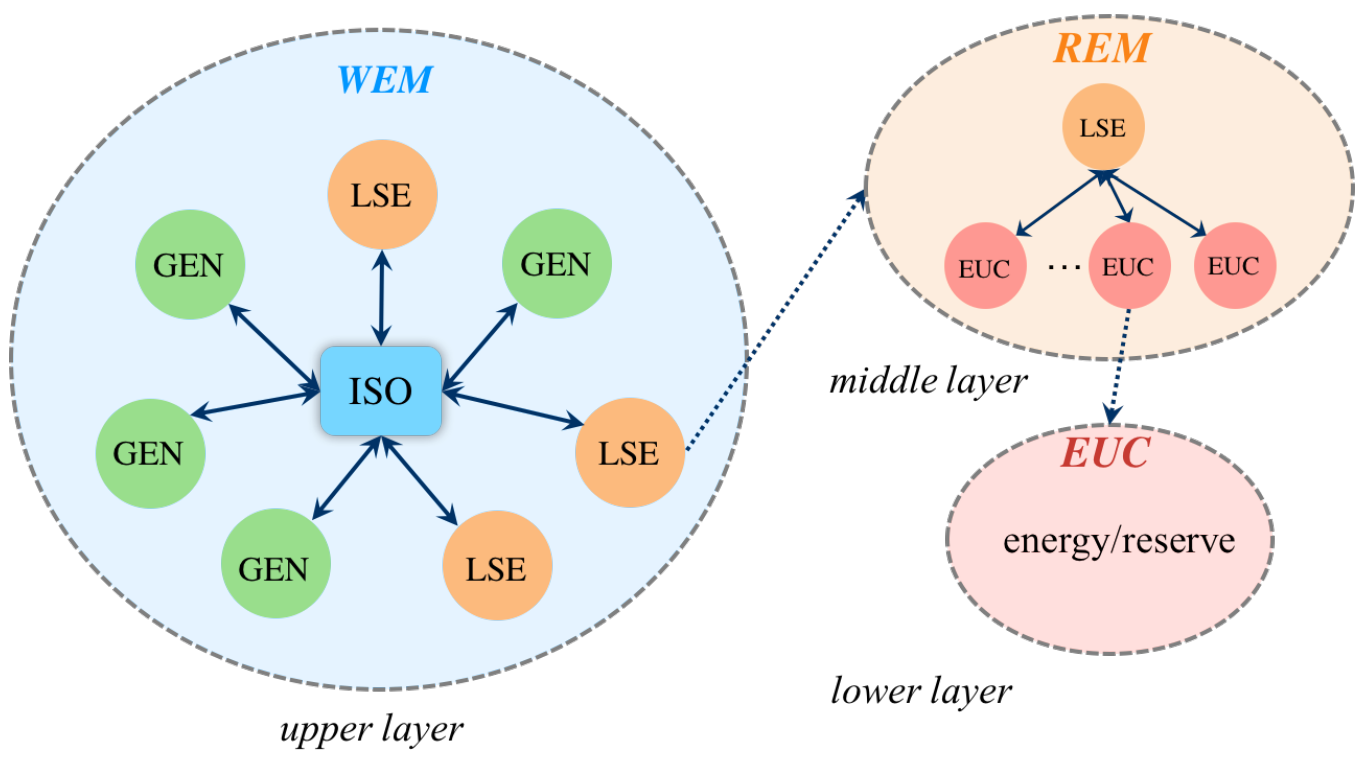

Fig. 1. Tri-layer structure of the electricity markets.

Consider a WEM where energy and spinning reserve are jointly cleared. A set of generators (GENs) submit energy and reserve offers in the WEM and a set of LSEs submit energy bids and reserve offers. Without loss of generality, we assume each energy offer/bid is associated with one reserve offer. For the purpose of simplicity, we refer to the energy and reserve offers from GENs as GEN offers, and energy bids and reserve offers from LSEs as LSE bids. Let $\mathcal{G} = \{1, \cdots, G\}$ be the index set of GEN offers and $\mathcal{D} = \{1, \cdots, D\}$ be the index set of LSE bids.

An offer/bid is defined to be a five-tuple that consists of marginal price, minimum quantity, and maximum quantity, for energy and reserve. Denote GEN offer $i$, $i \in \mathcal{G}$, as $(a_i^g, \underline{p}^g, \overline{p}^g, b_i^g, \overline{r}_i^g)$, where $a_i^g$ and $b_i^g$ are the marginal offer prices for energy and reserve, respectively, $\underline{p}_i^g$, $\overline{p}_i^g$ are the respective minimum and maximum active power supplies, and $\overline{r}_i^g$ is the maximum reserve. Define $\boldsymbol{a}^g = [a_1^g, \cdots, a_G^g]^\top$, $\underline{\boldsymbol{p}}^g = [\underline{p}_1^g, \cdots, \underline{p}_G^g]^\top$, $\overline{\boldsymbol{p}}^g = [\overline{p}_1^g, \cdots, \overline{p}_G^g]^\top$, $\boldsymbol{b}^g = [b_1^g, \cdots, b_G^g]^\top$, and $\overline{\boldsymbol{r}}^g = [\overline{r}_1^g, \cdots, \overline{r}_G^g]^\top$. Similarly, denote LSE bid $i$, $i \in \mathcal{D}$, as $(a_i^d, \underline{p}^d, \overline{p}^d, b_i^d, \overline{r}_i^d)$, where $a_i^d$ and $b_i^d$ are the marginal bid price for energy and the marginal offer price for reserve, respectively, $\underline{p}_i^d$, $\overline{p}_i^d$ are the respective minimum and maximum power demands, and $\overline{r}_i^d$ is the maximum reserve. Define $\boldsymbol{a}^d = [a_1^d, \cdots, a_D^d]^\top$, $\underline{\boldsymbol{p}}^d = [\underline{p}_1^d, \cdots, \underline{p}_D^d]^\top$, $\overline{\boldsymbol{p}}^d = [\overline{p}_1^d, \cdots, \overline{p}_D^d]^\top$, $\boldsymbol{b}^d = [b_1^d, \cdots, b_D^d]^\top$, $\overline{\boldsymbol{r}}^d = [\overline{r}_1^d, \cdots, \overline{r}_D^d]^\top$. Without loss of generality, we assume the minimum reserve is zero.

Let $\boldsymbol{\Phi}^g \in \mathbb{R}^{N \times G}$ be the GEN offer to bus mapping matrix, the $(i, j)^{th}$ entry of which is 1 if the $j^{th}$ offer from GEN is associated with bus $i$. The LSE bid to bus mapping matrix, denoted by $\boldsymbol{\Phi}^d \in \mathbb{R}^{N \times D}$, is defined in a similar way. Note that there is one and only one non-zero entry in each column, i.e., each offer/bid is only allowed to be submitted at one bus.

The WEM clearing problem for a single period,[1] which determines the cleared quantities, $\boldsymbol{p}^g$, $\boldsymbol{p}^d$, $\boldsymbol{r}^g$, $\boldsymbol{r}^d$, can be formulated as follows:

$$\underset{\boldsymbol{p}^g, \boldsymbol{p}^d, \boldsymbol{r}^g, \boldsymbol{r}^d}{\text{maximize}} \quad -(\boldsymbol{a}^g)^\top \boldsymbol{p}^g + (\boldsymbol{a}^d)^\top \boldsymbol{p}^d - (\boldsymbol{b}^g)^\top \boldsymbol{r}^g - (\boldsymbol{b}^d)^\top \boldsymbol{r}^d$$

subject to

$$(\boldsymbol{1} - \boldsymbol{\Lambda})^\top (\boldsymbol{\Phi}^g \boldsymbol{p}^g - \boldsymbol{\Phi}^d \boldsymbol{p}^d) = 0, \quad \leftrightarrow \quad \lambda \tag{1a}$$

$$-\overline{\boldsymbol{f}} \leq \boldsymbol{\Psi} (\boldsymbol{\Phi}^g \boldsymbol{p}^g - \boldsymbol{\Phi}^d \boldsymbol{p}^d) \leq \overline{\boldsymbol{f}}, \quad \leftrightarrow \quad \underline{\boldsymbol{\mu}}, \overline{\boldsymbol{\mu}} \tag{1b}$$

$$\boldsymbol{1}^\top \boldsymbol{r}^g + \boldsymbol{1}^\top \boldsymbol{r}^d \geq \underline{r}, \quad \leftrightarrow \quad \nu \tag{1c}$$

$$\underline{\boldsymbol{p}}^g \leq \boldsymbol{p}^g \leq \overline{\boldsymbol{p}}^g - \boldsymbol{r}^g, \quad \leftrightarrow \quad \underline{\boldsymbol{\rho}}^g, \overline{\boldsymbol{\rho}}^g \tag{1d}$$

$$\underline{\boldsymbol{p}}^d + \boldsymbol{r}^d \leq \boldsymbol{p}^d \leq \overline{\boldsymbol{p}}^d, \quad \leftrightarrow \quad \underline{\boldsymbol{\rho}}^d, \overline{\boldsymbol{\rho}}^d \tag{1e}$$

$$\boldsymbol{0} \leq \boldsymbol{r}^g \leq \overline{\boldsymbol{r}}^g, \quad \leftrightarrow \quad \underline{\boldsymbol{\eta}}^g, \overline{\boldsymbol{\eta}}^g \tag{1f}$$

$$\boldsymbol{0} \leq \boldsymbol{r}^d \leq \overline{\boldsymbol{r}}^d, \quad \leftrightarrow \quad \underline{\boldsymbol{\eta}}^d, \overline{\boldsymbol{\eta}}^d \tag{1g}$$

where $\underline{r}$ is the minimum reserve required, $\lambda$, $\underline{\boldsymbol{\mu}} \in \mathbb{R}^L$, $\overline{\boldsymbol{\mu}} \in \mathbb{R}^L$, $\nu$, $\underline{\boldsymbol{\rho}}^g \in \mathbb{R}^G$, $\overline{\boldsymbol{\rho}}^g \in \mathbb{R}^G$, $\underline{\boldsymbol{\rho}}^d \in \mathbb{R}^D$, $\overline{\boldsymbol{\rho}}^d \in \mathbb{R}^D$, $\underline{\boldsymbol{\eta}}^g \in \mathbb{R}^G$, $\overline{\boldsymbol{\eta}}^g \in \mathbb{R}^G$, $\underline{\boldsymbol{\eta}}^d \in \mathbb{R}^D$, and $\overline{\boldsymbol{\eta}}^d \in \mathbb{R}^D$ are dual variables associated with corresponding constraints, and $\boldsymbol{1}$ is the all-ones vector with an appropriate dimension. The objective function of the WEM clearing problem is the social welfare and the market clearing process aims to maximize the social welfare. Constraint (1a) is the power balance constraint with the explicit consideration of losses. Constraint (1b) is the line active power flow constraint. Constraint (1c) sets the minimum reserve requirements. Constraints (1d) to (1f) are the capacity constraints.

The WEM clearing problem defined by (1) determines the clearing prices for energy and reserve, as well as the quantities cleared from each offer and bid. Specifically, the clearing price of the reserves is $\nu$, and the locational marginal price (LMP) of energy at bus $i \in \mathcal{N}$, denoted by $\pi_i$, is calculated using the dual variables as follows [21]:

$$\pi_i = \lambda - \Lambda_i \lambda + \sum_{j=1}^{L} (\underline{\mu}_j - \overline{\mu}_j) \Psi_{ji}. \tag{2}$$

### B. Middle Layer: Retail Electricity Market

Consider a LSE that submits a set of energy bids and reserve offers denoted by $\tilde{\mathcal{D}} \subseteq \mathcal{D}$ with $|\tilde{\mathcal{D}}| = \tilde{D}$, where $|\cdot|$ denotes the cardinality of a set. Assume the LSE submits bids at a set of buses indexed by $\tilde{\mathcal{N}} = \{i_1, \cdots, i_n\} \subseteq \mathcal{N}$. Partition $\tilde{\mathcal{D}}$ into $n$ non-overlapping subsets, $\tilde{\mathcal{D}}_1, \cdots, \tilde{\mathcal{D}}_n$, such that $\tilde{\mathcal{D}} = \tilde{\mathcal{D}}_1 \cup \cdots \cup \tilde{\mathcal{D}}_n$ and bids in $\tilde{\mathcal{D}}_k$ are submitted at bus $i_k$, where the LMP is $\pi_{i_k}$. If a bid of the LSE gets cleared in the WEM, it will be charged at LMP for energy and paid at $\nu$ for reserve. From the LSE's perspective, $\pi_{i_k}$, for any $i_k \in \tilde{\mathcal{N}}$, $\nu$, $p_j^d$ and $r_j^d$ for any $j \in \tilde{\mathcal{D}}$ are functions of its bids defined through

[1]The single-period market clearing problem can be extended to a multi-period market clearing problem simply by including the same constraints for each single period and additionally imposing ramping constrains.

(1). The LSE may be able to change its clearing results by adjusting its bids. We assume the minimum and maximum quantities are fixed and the only decision variables in a bid are the marginal prices.

Assume the LSE operates a REM with a set of EUCs indexed by $\mathcal{C} = \{1, \cdots, C\}$. Partition $\mathcal{C}$ into $n$ non-overlapping subsets, $\mathcal{C}_1, \cdots, \mathcal{C}_n$, such that $\mathcal{C} = \mathcal{C}_1 \cup \cdots \cup \mathcal{C}_n$ and EUCs in $\mathcal{C}_k$ are associated with bus $i_k$ in the transmission network. We refer to $\mathcal{C}_k$ as pricing area $k$. Define $C_k = |\mathcal{C}_k|$. The LSE needs to determine the energy price $\alpha_k$ and the reserve price $\beta_k$ for EUCs in $\mathcal{C}_k$, for $k = 1, \ldots, n$. For a profit-seeking LSE, $\alpha_k$ and $\beta_k$ are chosen in such a way that its revenue is maximized. The price signals are sent to each EUC, which will then adjust its energy consumption and the reserve provision. Let $p_i^{(k)}$ and $r_i^{(k)}$ denote the respective energy demand and reserve provision from the $i^{th}$ EUC in $\mathcal{C}_k$. Define $\boldsymbol{p}^{(k)} = [p_1^{(k)}, \cdots, p_{C^k}^{(k)}]^\top$ and $\boldsymbol{r}^{(k)} = [r_1^{(k)}, \cdots, r_{C^k}^{(k)}]^\top$. From the LSE's perspective, $\boldsymbol{p}^{(k)}$ and $\boldsymbol{r}^{(k)}$ are functions of $\alpha_k$ and $\beta_k$ defined through an optimization problem that we will discuss in more details in next section.

The LSE has to determine $\alpha_k$ and $\beta_k$ for $k = 1 \cdots, n$, as well as $a_j^d$ and $b_j^d$ for all $j \in \tilde{\mathcal{D}}$, so as to maximize its total revenue from energy and reserve while ensuring that the energy purchased from the WEM and that sold in the REM balance at each bus in the transmission network, so do the reserve sold in the WEM and that purchased in the REM. This problem can be formulated as the following bi-level programming problem:

$$\underset{\substack{a_j^d, b_j^d, \forall j \in \tilde{\mathcal{D}} \\ \alpha_k, \beta_k, k=1\cdots,n}}{\text{maximize}} \quad \sum_{k=1}^{n} \sum_{j \in \tilde{\mathcal{D}}_k} \left( (\alpha_k - \pi_{i_k}) p_j^d + (\nu - \beta_k) r_j^d \right)$$

subject to

$$\pi_{i_k} = \lambda - \Lambda_{i_k}\lambda + \sum_{j=1}^{L} (\underline{\mu}_j - \overline{\mu}_j)\Psi_{j i_k}, \forall k, \tag{3a}$$

$$\sum_{j \in \tilde{\mathcal{D}}_k} p_j^d = \mathbf{1}^\top \boldsymbol{p}^{(k)}, \forall k, \tag{3b}$$

$$\sum_{j \in \tilde{\mathcal{D}}_k} r_j^d = \mathbf{1}^\top \boldsymbol{r}^{(k)}, \forall k, \tag{3c}$$

$$\underline{\alpha} \le \alpha_k \le \overline{\alpha}, \ \ \overline{\beta} \le \beta_k \le \overline{\beta}, \ \forall k, \tag{3d}$$

$$(1), \ (6), \ \forall k. \tag{3e}$$

where $\underline{\alpha}$, $\overline{\alpha}$, $\underline{\beta}$, $\overline{\beta}$ are price caps. The first term in the objective function is the total profit from energy trading and the second term is the total profit from reserve trading. Constraint (3a) is the expression to calculate the LMP. Constraints (3b) and (3c) ensure the respective energy and reserve balances. Constraint (3d) restricts the energy and reserve prices within a pre-specified range. The impacts from the WEM and the EUCs are represented respectively by (1) and (6) (to be detailed in next section) in constraint (3e). Specifically, $\lambda$, $\underline{\mu}_j$, and $\overline{\mu}_j$ are determined through (1), while $\boldsymbol{p}^{(k)}$ and $\boldsymbol{r}^{(k)}$ are determined through (6).

### C. Lower Layer: End User Customers

Each EUC purchases energy and sells reserve in the REM. The objective of the $i^{\text{th}}$ EUC in $\mathcal{C}_k$ is to maximize its own benefit given the prices $\alpha_k$ and $\beta_k$. The EUC determines its energy demand $p_i^{(k)}$ and reserve provision $r_i^{(k)}$ by solving the following problem:

$$\underset{p_i^{(k)}, r_i^{(k)}}{\text{maximize}} \quad g_i^{(k)}(p_i^{(k)}) - \alpha_k p_i^{(k)} + \beta_k r_i^{(k)} - h_i^{(k)}(r_i^{(k)})$$

subject to

$$\underline{p}_i^{(k)} + r_i^{(k)} \le p_i^{(k)} \le \overline{p}_i^{(k)}, \tag{4a}$$

$$0 \le r_i^{(k)} \le \overline{r}_i^{(k)}, \tag{4b}$$

where $g_i^{(k)}$ is the energy benefit function of the $i^{\text{th}}$ EUC in $\mathcal{C}_k$, and $h_i^{(k)}$ is the reserve cost function. To simplify the EUC modeling, we assume $g_i^{(k)}$ and $h_i^{(k)}$ are piecewise linear functions, i.e., $g_i^{(k)}(p_i^{(k)}) = \sum_{j=1}^{n_i^{(k)}} c_{i,j}^{(k)} x_{i,j}^{(k)}$ and $h_i^{(k)}(r_i^{(k)}) = \sum_{j=1}^{n_i^{(k)}} d_{i,j}^{(k)} y_{i,j}^{(k)}$, where $x_{i,j}^{(k)} \in [\underline{x}_{i,j}^{(k)}, \overline{x}_{i,j}^{(k)}]$ and $y_{i,j}^{(k)} \in [0, \overline{y}_{i,j}^{(k)}]$ are the $j^{\text{th}}$ block of energy demand and reserve provision of the $i^{\text{th}}$ EUC in $\mathcal{C}_k$, respectively, and $n_i^{(k)}$ is the number of blocks. Note that $p_i^{(k)} = \sum_{j=1}^{n_j} x_{i,j}^{(k)}$, $r_i^{(k)} = \sum_{j=1}^{n_j} y_{i,j}^{(k)}$, $\sum_{j=1}^{n_i^{(k)}} \underline{x}_{i,j}^{(k)} = \underline{p}_i^{(k)}$, $\sum_{j=1}^{n_i^{(k)}} \overline{x}_{i,j}^{(k)} = \overline{p}_i^{(k)}$, and $\sum_{j=1}^{n_i^{(k)}} \overline{y}_{i,j}^{(k)} = \overline{r}_i^{(k)}$. Then, (4) becomes

$$\underset{x_{i,j}^{(k)}, y_{i,j}^{(k)}}{\text{maximize}} \quad \sum_{j=1}^{n_i^{(k)}} \left( (c_{i,j}^{(k)} - \alpha_k) x_{i,j}^{(k)} + (\beta_k - d_{i,j}^{(k)}) y_{i,j}^{(k)} \right)$$

subject to

$$\underline{x}_{i,j}^{(k)} + y_{i,j}^{(k)} \le x_{i,j}^{(k)} \le \overline{x}_{i,j}^{(k)}, j = 1, \cdots, n_j, \tag{5a}$$

$$0 \le y_{i,j}^{(k)} \le \overline{y}_{i,j}^{(k)}, j = 1, \cdots, n_j. \tag{5b}$$

Problem (5) solved by each EUC is independent of problems solved by others. The set of problem of EUCs in $\mathcal{C}_k$ can be aggregately written in the following vector-form:

$$\underset{\boldsymbol{x}^{(k)}, \boldsymbol{y}^{(k)}}{\text{maximize}} \quad (\boldsymbol{c}^{(k)})^\top \boldsymbol{x}^{(k)} - \alpha_k \mathbf{1}^\top \boldsymbol{x}^{(k)} + \beta_k \mathbf{1}^\top \boldsymbol{y}^{(k)} - (\boldsymbol{d}^{(k)})^\top \boldsymbol{y}^{(k)}$$

subject to

$$\underline{\boldsymbol{x}}^{(k)} + \boldsymbol{y}^{(k)} \le \boldsymbol{x}^{(k)} \le \overline{\boldsymbol{x}}^{(k)}, \ \leftrightarrow \ \underline{\boldsymbol{\gamma}}^{(k)}, \overline{\boldsymbol{\gamma}}^{(k)} \tag{6a}$$

$$\mathbf{0} \le \boldsymbol{y}^{(k)} \le \overline{\boldsymbol{y}}^{(k)}, \ \leftrightarrow \ \underline{\boldsymbol{\zeta}}^{(k)}, \overline{\boldsymbol{\zeta}}^{(k)} \tag{6b}$$

where $\boldsymbol{c}^{(k)}$ and $\boldsymbol{d}^{(k)}$ are the coefficients of the energy benefit functions and reserve cost functions, $\underline{\boldsymbol{\gamma}}^{(k)}$, $\overline{\boldsymbol{\gamma}}^{(k)}$, $\underline{\boldsymbol{\zeta}}^{(k)}$ and $\overline{\boldsymbol{\zeta}}^{(k)}$ are dual variables associated with corresponding constraints.

### D. Discussion

It is in general inconvenient for an LSE to know the exact utility functions of the EUCs. However, an LSE could assume some form of the utility functions and further approximate them using piecewise linear functions, such as the method in the prior work [11]. The parameters of such functions may be estimated from actual data using some forecast techniques for energy consumption, e.g., in [22]. Yet, the estimation of the parameters in these functions is another interesting problem, which is worth of being investigated in future works.

Another interesting topic is how to determine the constraints when an EUC makes its decision in reality. As a basic proof-of-concept for the tri-layer market structure, we use the common box constraints of power consumption, which is based on an assumption that the utility function of an EUC can be estimated by the LSE using collected data, and then these box constraints can be well captured by the piecewise linear utility function. For EUCs with specific constraints, we believe a more complete and detailed model for the EUCs is indeed necessary and deserve further investigation in the future.

## III. Joint Bidding and Pricing Problem Formulation

In this section, we formulate the joint bidding and pricing problem of the LSE into a single-level mixed integer programming (MIP) problem by plugging the KKT conditions of the upper layer problem (1) and the lower layer problem (6) into the middle layer problem (3) .

### A. Characterization of WEM Clearing Problem

The clearing results of the WEM can be characterized by the KKT conditions of problem (1) derived as follows:

*Stationarity*:

$$\begin{aligned} \boldsymbol{a}^g - (\boldsymbol{\Phi}^g)^\top(\lambda(\mathbf{1}-\boldsymbol{\Lambda}) + \boldsymbol{\Psi}^\top(\underline{\boldsymbol{\mu}} - \overline{\boldsymbol{\mu}})) - \underline{\boldsymbol{\rho}}^g + \overline{\boldsymbol{\rho}}^g = \mathbf{0}, \\ -\boldsymbol{a}^d + (\boldsymbol{\Phi}^d)^\top(\lambda(\mathbf{1}-\boldsymbol{\Lambda}) + \boldsymbol{\Psi}^\top(\underline{\boldsymbol{\mu}} - \overline{\boldsymbol{\mu}})) - \underline{\boldsymbol{\rho}}^d + \overline{\boldsymbol{\rho}}^d = \mathbf{0}, \\ \boldsymbol{b}^g - \nu\mathbf{1} + \overline{\boldsymbol{\rho}}^g - \underline{\boldsymbol{\eta}}^g + \overline{\boldsymbol{\eta}}^g = \mathbf{0}, \\ \boldsymbol{b}^d - \nu\mathbf{1} + \underline{\boldsymbol{\rho}}^d - \underline{\boldsymbol{\eta}}^d + \overline{\boldsymbol{\eta}}^d = \mathbf{0}. \end{aligned} \tag{7}$$

*Primal feasibility*:

$$\text{(1a)} - \text{(1g)}. \tag{8}$$

*Dual feasibility*:

$$\begin{aligned} \underline{\boldsymbol{\mu}} \geq \mathbf{0}, \overline{\boldsymbol{\mu}} \geq \mathbf{0}, \nu \geq 0, \underline{\boldsymbol{\rho}}^g \geq \mathbf{0}, \overline{\boldsymbol{\rho}}^g \geq \mathbf{0}, \\ \underline{\boldsymbol{\rho}}^d \geq \mathbf{0}, \overline{\boldsymbol{\rho}}^d \geq \mathbf{0}, \overline{\boldsymbol{\eta}}^g \geq \mathbf{0}, \overline{\boldsymbol{\eta}}^g \geq \mathbf{0}, \underline{\boldsymbol{\eta}}^d \geq \mathbf{0}, \overline{\boldsymbol{\eta}}^d \geq \mathbf{0}. \end{aligned} \tag{9}$$

*Complementary slackness*:

$$\begin{aligned} \underline{\boldsymbol{\mu}} \circ (\boldsymbol{\Psi}(\boldsymbol{\Phi}^g \boldsymbol{p}^g - \boldsymbol{\Phi}^d \boldsymbol{p}^d) + \overline{\boldsymbol{f}}) = \mathbf{0}, \\ \overline{\boldsymbol{\mu}} \circ (\boldsymbol{\Psi}(\boldsymbol{\Phi}^g \boldsymbol{p}^g - \boldsymbol{\Phi}^d \boldsymbol{p}^d) - \overline{\boldsymbol{f}}) = \mathbf{0}, \\ \nu(\mathbf{1}^\top \boldsymbol{r}^g + \mathbf{1}^\top \boldsymbol{r}^d - \underline{r}) = 0, \quad \underline{\boldsymbol{\rho}}^g \circ (\boldsymbol{p}^g - \underline{\boldsymbol{p}}^g) = \mathbf{0}, \\ \overline{\boldsymbol{\rho}}^g \circ (\boldsymbol{p}^g + \boldsymbol{r}^g - \overline{\boldsymbol{p}}^g) = \mathbf{0}, \quad \overline{\boldsymbol{\eta}}^g \circ (\boldsymbol{r}^g - \overline{\boldsymbol{r}}^g) = \mathbf{0}, \\ \underline{\boldsymbol{\rho}}^d \circ (\boldsymbol{p}^d - \boldsymbol{r}^d - \underline{\boldsymbol{p}}^d) = \mathbf{0}, \quad \overline{\boldsymbol{\rho}}^d \circ (\boldsymbol{p}^d - \overline{\boldsymbol{p}}^d) = \mathbf{0}, \\ \overline{\boldsymbol{\eta}}^d \circ (\boldsymbol{r}^d - \overline{\boldsymbol{r}}^d) = \mathbf{0}, \quad \underline{\boldsymbol{\eta}}^g \circ \boldsymbol{r}^g = \mathbf{0}, \quad \underline{\boldsymbol{\eta}}^d \circ \boldsymbol{r}^d = \mathbf{0}, \end{aligned} \tag{10}$$

where $\circ$ indicates element-wise product, and $\mathbf{0}$ is the all-zeros vector with an appropriate dimension. Note that constraints (7) to (9) are linear. The only non-linear constraints are the complementary slackness conditions, which are indeed bilinear with respect to the unknown primal and dual variables.

### B. Characterization of EUC Benefit Maximization Problem

The benefit maximization problem of EUCs in $\mathcal{C}_k$ can be completely characterized by the KKT conditions of problem (6) derived as follows:

*Stationarity*:

$$\begin{aligned} -\boldsymbol{c}^{(k)} + \alpha_k \mathbf{1} - \underline{\boldsymbol{\gamma}}^{(k)} + \overline{\boldsymbol{\gamma}}^{(k)} = \mathbf{0}, \\ \boldsymbol{d}^{(k)} - \beta_k \mathbf{1} + \underline{\boldsymbol{\gamma}}^{(k)} - \underline{\boldsymbol{\zeta}}^{(k)} + \overline{\boldsymbol{\zeta}}^{(k)} = \mathbf{0}. \end{aligned} \tag{11}$$

*Primal feasibility*:

$$\text{(6a)} - \text{(6b)}. \tag{12}$$

*Dual feasibility*:

$$\underline{\boldsymbol{\gamma}}^{(k)} \geq \mathbf{0}, \overline{\boldsymbol{\gamma}}^{(k)} \geq \mathbf{0}, \underline{\boldsymbol{\zeta}}^{(k)} \geq \mathbf{0}, \overline{\boldsymbol{\zeta}}^{(k)} \geq \mathbf{0}. \tag{13}$$

*Complementary slackness*:

$$\begin{aligned} \underline{\boldsymbol{\gamma}}^{(k)} \circ (\boldsymbol{x}^{(k)} - \boldsymbol{y}^{(k)} - \underline{\boldsymbol{x}}^{(k)}) = \mathbf{0}, \\ \overline{\boldsymbol{\gamma}}^{(k)} \circ (\boldsymbol{x}^{(k)} - \overline{\boldsymbol{x}}^{(k)}) = \mathbf{0}, \\ \underline{\boldsymbol{\zeta}}^{(k)} \circ \boldsymbol{y}^{(k)} = \mathbf{0}, \quad \overline{\boldsymbol{\zeta}}^{(k)} \circ (\boldsymbol{y}^{(k)} - \overline{\boldsymbol{y}}^{(k)}) = \mathbf{0}. \end{aligned} \tag{14}$$

Similar to the stationarity condition in the WEM clearing problem, equation (11) is linear with respect to marginal prices of benefit/cost and the dual variables. The only non-linear constraints are the complementary slackness conditions, which are bilinear with respect to the primal and dual variables.

### C. Formulation of LSE Joint Bidding and Pricing Problem

Replacing (3e) in (3) by the KKT conditions III-A and III-B, we can obtain the joint bidding and pricing problem of LSE $i$ as follows:

$$\underset{\substack{a_j^d, b_j^d, \forall j \in \tilde{\mathcal{D}} \\ \alpha_k, \beta_k, k=1\cdots,n}}{\text{maximize}} \quad \sum_{k=1}^{n} \sum_{j \in \tilde{\mathcal{D}}_k} ((\alpha_k - \pi_{i_k}) p_j^d + (\nu - \beta_k) r_j^d)$$

subject to

$$\text{(3a)} - \text{(3d)}, \ \text{(7)} - \text{(14)}, \tag{15}$$

which has a bilinear constraints and objective function, with continuous/binary decision variables.

The formulation of the LSE joint bidding and pricing problem requires the LSE hypothesizes offer/bid functions of other WEM participants, and benefit and cost functions of EUCs that participate in the REM as in existing works [10], [19]. While the estimation or learning of offer/bid functions is itself an interesting and valuable problem, it is out the scope of this paper. We refer the interested readers to [23] and references therein for more details.

## IV. LSE Joint Bidding and Pricing Problem Linearization

The MIP in Section II has bilinear constraints and objective function, which poses significant challenge to solve it. In this section, we linearize the bilinear terms and transform the complete pricing problem (15) into an MILP problem.

### A. Linearization of Complementary Slackness Conditions

Using the techniques from [24], the bilinear complementary slackness conditions (III-A) can be linearized as follows:

$$\begin{aligned}
&\boldsymbol{\Psi}(\boldsymbol{\Phi}^g\boldsymbol{p}^g-\boldsymbol{\Phi}^d\boldsymbol{p}^d)+\overline{\boldsymbol{f}}\le M(\mathbf{1}-\boldsymbol{z}^{\underline{\mu}}),\ \underline{\boldsymbol{\mu}}\le M\boldsymbol{z}^{\underline{\mu}},\\
&\overline{\boldsymbol{f}}-\boldsymbol{\Psi}(\boldsymbol{\Phi}^g\boldsymbol{p}^g-\boldsymbol{\Phi}^d\boldsymbol{p}^d)\le M(\mathbf{1}-\boldsymbol{z}^{\overline{\mu}}),\ \overline{\boldsymbol{\mu}}\le M\boldsymbol{z}^{\overline{\mu}},\\
&\underline{r}-\mathbf{1}^\top\boldsymbol{r}^g+\mathbf{1}^\top\boldsymbol{r}^d\le M(1-z^\nu),\ \nu\le Mz^\nu,\\
&\boldsymbol{p}^g-\underline{\boldsymbol{p}}^g\le M(\mathbf{1}-\boldsymbol{z}^{\underline{\rho}^g}),\ \underline{\boldsymbol{\rho}}^g\le M\boldsymbol{z}^{\underline{\rho}^g},\\
&\overline{\boldsymbol{p}}^g-\boldsymbol{p}^g-\boldsymbol{r}^g\le M(\mathbf{1}-\boldsymbol{z}^{\overline{\rho}^g}),\ \overline{\boldsymbol{\rho}}^g\le M\boldsymbol{z}^{\overline{\rho}^g},\\
&\boldsymbol{p}^d-\boldsymbol{r}^d-\underline{\boldsymbol{p}}^d\le M(\mathbf{1}-\boldsymbol{z}^{\underline{\rho}^d}),\ \underline{\boldsymbol{\rho}}^d\le M\boldsymbol{z}^{\underline{\rho}^d},\\
&\overline{\boldsymbol{p}}^d-\boldsymbol{p}^d\le M(\mathbf{1}-\boldsymbol{z}^{\overline{\rho}^d}),\ \overline{\boldsymbol{\rho}}^d\le M\boldsymbol{z}^{\overline{\rho}^d},\\
&\boldsymbol{r}^g\le M(\mathbf{1}-\boldsymbol{z}^{\underline{\eta}^g}),\ \underline{\boldsymbol{\eta}}^g\le M\boldsymbol{z}^{\underline{\eta}^g},\\
&\overline{\boldsymbol{r}}^g-\boldsymbol{r}^g\le M(\mathbf{1}-\boldsymbol{z}^{\overline{\eta}^g}),\ \overline{\boldsymbol{\eta}}^g\le M\boldsymbol{z}^{\overline{\eta}^g},\\
&\boldsymbol{r}^d\le M(\mathbf{1}-\boldsymbol{z}^{\underline{\eta}^d}),\ \underline{\boldsymbol{\eta}}^d\le M\boldsymbol{z}^{\underline{\eta}^d},\\
&\overline{\boldsymbol{r}}^d-\boldsymbol{r}^d\le M(\mathbf{1}-\boldsymbol{z}^{\overline{\eta}^d}),\ \overline{\boldsymbol{\eta}}^d\le M\boldsymbol{z}^{\overline{\eta}^d},\\
&\left.\begin{matrix}\boldsymbol{z}^{\underline{\mu}},\boldsymbol{z}^{\overline{\mu}},z^{\mu},\boldsymbol{z}^{\underline{\rho}^g},\boldsymbol{z}^{\overline{\rho}^g},\\ \boldsymbol{z}^{\underline{\rho}^d},\boldsymbol{z}^{\overline{\rho}^d},\boldsymbol{z}^{\underline{\eta}^g},\boldsymbol{z}^{\overline{\eta}^g},\boldsymbol{z}^{\underline{\eta}^d},\boldsymbol{z}^{\overline{\eta}^d}\end{matrix}\right.\in\{0,1\}.
\end{aligned}\tag{16}$$

where $M$ is a big number [2].

Similarly, the bilinear complementary slackness conditions (14) can be linearized as follows:

$$\begin{aligned}
&\boldsymbol{x}^{(k)}-\boldsymbol{y}^{(k)}-\underline{\boldsymbol{x}}^{(k)}\le M(\mathbf{1}-\boldsymbol{z}^{\underline{\gamma}^{(k)}}),\ \underline{\boldsymbol{\gamma}}^{(k)}\le M\boldsymbol{z}^{\underline{\gamma}^{(k)}},\\
&\overline{\boldsymbol{x}}^{(k)}-\boldsymbol{x}^{(k)}\le M(\mathbf{1}-\boldsymbol{z}^{\overline{\gamma}^{(k)}}),\ \overline{\boldsymbol{\gamma}}^{(k)}\le M\boldsymbol{z}^{\overline{\gamma}^{(k)}},\\
&\boldsymbol{y}^{(k)}\le M(\mathbf{1}-\boldsymbol{z}^{\underline{\zeta}^{(k)}}),\ \underline{\boldsymbol{\zeta}}^{(k)}\le M\boldsymbol{z}^{\underline{\zeta}^{(k)}},\\
&\overline{\boldsymbol{y}}^{(k)}-\boldsymbol{y}^{(k)}\le M(\mathbf{1}-\boldsymbol{z}^{\overline{\zeta}^{(k)}}),\ \overline{\boldsymbol{\zeta}}^{(k)}\le M\boldsymbol{z}^{\overline{\zeta}^{(k)}},\\
&\boldsymbol{z}^{\underline{\gamma}^{(k)}},\boldsymbol{z}^{\overline{\gamma}^{(k)}},\boldsymbol{z}^{\underline{\zeta}^{(k)}},\boldsymbol{z}^{\overline{\zeta}^{(k)}}\in\{0,1\}.
\end{aligned}\tag{17}$$

So far, all constraints in (15) become linear and the only nonlinear term left in (15) is the bilinear objective function, which is linearized in the next section.

### B. Linearization of Bilinear Objective Function

Since (6) is a linear programming problem, the strong duality holds, i.e.,

$$\begin{aligned}
&-(\boldsymbol{c}^{(k)})^\top\boldsymbol{x}^{(k)}+\alpha_k\mathbf{1}^\top\boldsymbol{x}^{(k)}-\beta_k\mathbf{1}^\top\boldsymbol{y}^{(k)}+(\boldsymbol{d}^{(k)})^\top\boldsymbol{y}^{(k)}\\
&=(\underline{\boldsymbol{\gamma}}^{(k)})^\top\underline{\boldsymbol{x}}^{(k)}-(\overline{\boldsymbol{\gamma}}^{(k)})^\top\overline{\boldsymbol{x}}^{(k)}-(\overline{\boldsymbol{\zeta}}^{(k)})^\top\overline{\boldsymbol{y}}^{(k)}.
\end{aligned}\tag{18}$$

As such, from (3b), (3c), and (18), the following relation holds:

$$\begin{aligned}
&\sum_{k=1}^n\sum_{j\in\tilde{\mathcal{D}}_k}(\alpha_kp_j^d-\beta_kr_j^d)=\sum_{k=1}^n((\boldsymbol{c}^{(k)})^\top\boldsymbol{x}^{(k)}-(\boldsymbol{d}^{(k)})^\top\boldsymbol{y}^{(k)}\\
&\quad+(\underline{\boldsymbol{\gamma}}^{(k)})^\top\underline{\boldsymbol{x}}^{(k)}-(\overline{\boldsymbol{\gamma}}^{(k)})^\top\overline{\boldsymbol{x}}^{(k)}-(\overline{\boldsymbol{\zeta}}^{(k)})^\top\overline{\boldsymbol{y}}^{(k)}).
\end{aligned}\tag{19}$$

[2]Note that despite we use the same notation $M$ to represent a big number, the value of $M$ may be chosen differently in different constraints.

Similarly, the strong duality holds since (1) is a linear programming problem, i.e.,

$$\begin{aligned}
&(\boldsymbol{a}^g)^\top\boldsymbol{p}^g-(\boldsymbol{a}^d)^\top\boldsymbol{p}^d+(\boldsymbol{b}^g)^\top\boldsymbol{r}^g+(\boldsymbol{b}^d)^\top\boldsymbol{r}^d=\\
&\quad-(\underline{\boldsymbol{\mu}}+\overline{\boldsymbol{\mu}})^\top\overline{\boldsymbol{f}}+\nu\underline{r}+(\underline{\boldsymbol{\rho}}^g)^\top\underline{\boldsymbol{p}}^g-(\overline{\boldsymbol{\rho}}^g)^\top\overline{\boldsymbol{p}}^g\\
&\quad+(\underline{\boldsymbol{\rho}}^d)^\top\underline{\boldsymbol{p}}^d-(\overline{\boldsymbol{\rho}}^d)^\top\overline{\boldsymbol{p}}^d-(\overline{\boldsymbol{\eta}}^g)^\top\overline{\boldsymbol{r}}^g-(\overline{\boldsymbol{\eta}}^d)^\top\overline{\boldsymbol{r}}^d
\end{aligned}\tag{20}$$

which can be also written as:

$$\begin{aligned}
&\sum_{j\in\tilde{\mathcal{D}}}(-a_j^dp_j^d+b_j^dr_j^d)=\sum_{j\in\mathcal{D}-\tilde{\mathcal{D}}}(a_j^dp_j^d-b_j^dr_j^d)\\
&-(\boldsymbol{a}^g)^\top\boldsymbol{p}^g-(\boldsymbol{b}^g)^\top\boldsymbol{r}^g-(\underline{\boldsymbol{\mu}}+\overline{\boldsymbol{\mu}})^\top\overline{\boldsymbol{f}}+\nu\underline{r}+(\underline{\boldsymbol{\rho}}^g)^\top\underline{\boldsymbol{p}}^g\\
&-(\overline{\boldsymbol{\rho}}^g)^\top\overline{\boldsymbol{p}}^g+(\underline{\boldsymbol{\rho}}^d)^\top\underline{\boldsymbol{p}}^d-(\overline{\boldsymbol{\rho}}^d)^\top\overline{\boldsymbol{p}}^d-(\overline{\boldsymbol{\eta}}^g)^\top\overline{\boldsymbol{r}}^g-(\overline{\boldsymbol{\eta}}^d)^\top\overline{\boldsymbol{r}}^d.
\end{aligned}\tag{21}$$

The following equalities hold from (7), $\forall j\in\tilde{\mathcal{D}}_k$:

$$\pi_{i_k}=a_j^d+\underline{\rho}_j^d-\overline{\rho}_j^d,\quad \nu=b_i^d+\underline{\rho}_j^d-\underline{\eta}_j^d+\overline{\eta}_j^d.\tag{22}$$

The following equalities hold from (III-A):

$$\begin{aligned}
&\underline{\rho}_j^d(p_j^d-r_j^d)=\underline{\rho}_j^d\underline{p}_j^d,\quad \overline{\rho}_j^dp_j^d=\overline{\rho}_j^d\overline{p}_j^d,\\
&\underline{\eta}_j^dr_j^d=0,\quad \overline{\eta}_j^dr_j^d=\overline{\eta}_j^d\overline{r}_j^d.
\end{aligned}\tag{23}$$

From (21) – (23), the following equality holds:

$$\begin{aligned}
&\sum_{k=1}^n\sum_{j\in\tilde{\mathcal{D}}_k}(-\pi_{i_k}p_j^d+\nu r_j^d)=\sum_{j\in\tilde{\mathcal{D}}}(-\underline{\rho}_j^d\underline{p}_j^d+\overline{\rho}_j^d\overline{p}_j^d+\overline{\eta}_j^d\overline{r}_j^d)\\
&\quad+\sum_{j\in\mathcal{D}-\tilde{\mathcal{D}}}(a_j^dp_j^d-b_j^dr_j^d)-(\boldsymbol{a}^g)^\top\boldsymbol{p}^g-(\boldsymbol{b}^g)^\top\boldsymbol{r}^g\\
&\quad-(\underline{\boldsymbol{\mu}}+\overline{\boldsymbol{\mu}})^\top\overline{\boldsymbol{f}}+\nu\underline{r}+(\underline{\boldsymbol{\rho}}^g)^\top\underline{\boldsymbol{p}}^g-(\overline{\boldsymbol{\rho}}^g)^\top\overline{\boldsymbol{p}}^g\\
&\quad+(\underline{\boldsymbol{\rho}}^d)^\top\underline{\boldsymbol{p}}^d-(\overline{\boldsymbol{\rho}}^d)^\top\overline{\boldsymbol{p}}^d-(\overline{\boldsymbol{\eta}}^g)^\top\overline{\boldsymbol{r}}^g-(\overline{\boldsymbol{\eta}}^d)^\top\overline{\boldsymbol{r}}^d.
\end{aligned}\tag{24}$$

which expresses $\sum_{k=1}^n\sum_{j\in\tilde{\mathcal{D}}_k}(-\pi_{i_k}p_j^d+\nu r_j^d)$ as a linear combination of decision variables.

Combining (19) and (24), the objective function of (15) is now completely linearized, and together with constraints (3a) – (3c), (7) – (9), (16), (11) – (13), (17), the LSE pricing problem in (15) is transformed into an MILP problem and can be solved efficiently using existing solvers.

## V. Numerical Simulation

In this section, we present several numerical examples to illustrate the application of the proposed methodology using the IEEE 9-bus system and the IEEE 118-bus system test cases. All simulation results are obtained using Matlab 2016b on a Laptop with the CPU as 2.6 GHz Intel Core i7 and the memory as 16 GB. The solver for solving the MILP problem is Gurobi v7.5.2 [25].

### A. Small System Test Case

We use a modified version of the IEEE 9-bus test case [26] as an illustrative example. Buses $1-3$ represent generators and buses $4-9$ represent load buses. We increase the line active power capacity so that no congestion on the transmission lines occur. Moreover, the active power losses are assumed to be negligible, thus the LF vector $\boldsymbol{\Lambda}$ is $\mathbf{0}$.

TABLE I
BASE CASE CLEARING RESULTS WITH/W.O. JOINT BIDDING/PRICING

| | $\boldsymbol{\alpha}^*$ (\$/MW) | | | $\boldsymbol{a}^{d*}$ (\$/MW) | | |
|---|---|---|---|---|---|---|
| | $\mathcal{C}_1$ | $\mathcal{C}_2$ | $\mathcal{C}_3$ | $\mathcal{C}_1$ | $\mathcal{C}_2$ | $\mathcal{C}_3$ |
| with | 35.094 | 35.325 | 36.299 | 22.408 | 22.408 | 22.408 |
| w.o. | | 23.885 | | | 23.885 | |

| | $\boldsymbol{\beta}^*$ (\$/MW) | | | $\boldsymbol{b}^{d*}$ (\$/MW) | | |
|---|---|---|---|---|---|---|
| | $\mathcal{C}_1$ | $\mathcal{C}_2$ | $\mathcal{C}_3$ | $\mathcal{C}_1$ | $\mathcal{C}_2$ | $\mathcal{C}_3$ |
| with | 6.180 | 5.009 | 5.706 | 5.924 | 5.924 | 5.924 |
| w.o. | | 6.307 | | | 6.307 | |

| | LMP $\pi^*$ (\$/MW) | profit (\$) | welfare (\$) |
|---|---|---|---|
| with | 22.408 | 2595.422 | 2578.448 |
| w.o. | 23.885 | 0 | 2863.647 |

Consider an LSE that manages 3 different pricing areas, associated with EUC sets $\mathcal{C}_1, \mathcal{C}_2, \mathcal{C}_3$, corresponding to buses $6, 7,$ and $8$, respectively. The other load buses $4, 5,$ and $9$ correspond to pricing areas $\mathcal{C}_4, \mathcal{C}_5,$ and $\mathcal{C}_6$ controlled by other LSEs. Without loss of generality, we assume that each pricing area submits one offer/bid block to the WEM. We assume each pricing area consists of 10 EUCs, with each EUC submitting 3 bids, i.e., $C_k = 10$ and $n_i^{(k)} = 3$, $k = 1, 2, 3$. The values of $\boldsymbol{c}^{(k)}$ for $k = 1, 2, 3$ are independent and identical distributed (i.i.d.) samples drawn uniformly from $[34, 36]$, $[35, 37]$, and $[36, 38]$, respectively. The values of $\boldsymbol{d}^{(k)}$ are i.i.d. samples drawn uniformly from $[4, 6]$, $[5, 7]$, and $[6, 8]$, respectively. The upper bounds of EUCs' energy demands and reserve provisions in area $k$, $\bar{\boldsymbol{x}}^{(k)}$ and $\bar{\boldsymbol{y}}^{(k)}$, are i.i.d. samples drawn uniformly from $[0.85/C_k, 1.15/C_k] \cdot \bar{p}_k$ and $[0.85/C_k, 1.15/C_k] \cdot \bar{r}_k$, respectively. In addition, the offers of other LSEs, i.e., $a_j^d$ for $j = 4, 5, 6$, are estimated to be i.i.d. samples drawn uniformly from $[20, 22]$, $[21, 23]$, and $[22, 24]$, respectively, and their bids, i.e., $b_j^d$ for $j = 4, 5, 6$, are also i.i.d. samples drawn uniformly from $[4, 6]$, $[5, 7]$, and $[6, 8]$ as the pricing areas $\mathcal{C}_1, \mathcal{C}_2, \mathcal{C}_3$. The marginal offer prices for energy $a_j^g$ for $j = 1, 2, 3$ are sampled uniformly from $[20, 22]$, $[21, 23]$, and $[22, 24]$, respectively, for the three generators; and the marginal offer prices for reserve $b_j^g$ for $j = 1, 2, 3$ are sampled uniformly from $[3, 5]$, $[4, 6]$, and $[5, 7]$. Note that the offers provided by the generators are in general lower than the bids submitted by the pricing areas, such that most of the demands can be cleared. We also set the upper bounds for the energy price $\alpha_k$ and reserve price $\beta_k$, i.e., $\overline{\alpha}$ and $\overline{\beta}$, to be 100 and 50, respectively. The lower bounds $\underline{\alpha}$ and $\underline{\beta}$ are both zero.

### B. Base Case Clearing Results

The market clearing results of our small test case are shown in Table I. The clearing price is identical at all buses since the line congestion does not occur in the base case and the active power losses are ignored. It is seen that the energy and reserve prices vary across the areas $\mathcal{C}_k$, even though the WEM clearing price $\pi$ is identical at all buses. In fact, the LSE makes profits by adapting the prices to the EUCs' preference. Supposedly if the EUCs' willingness to buy is higher than the market clearing price, the LSE simply sets $\alpha_k$ to be the lowest value in $\mathbf{c}^{(k)}$, to clear as many demands from the EUCs as possible.

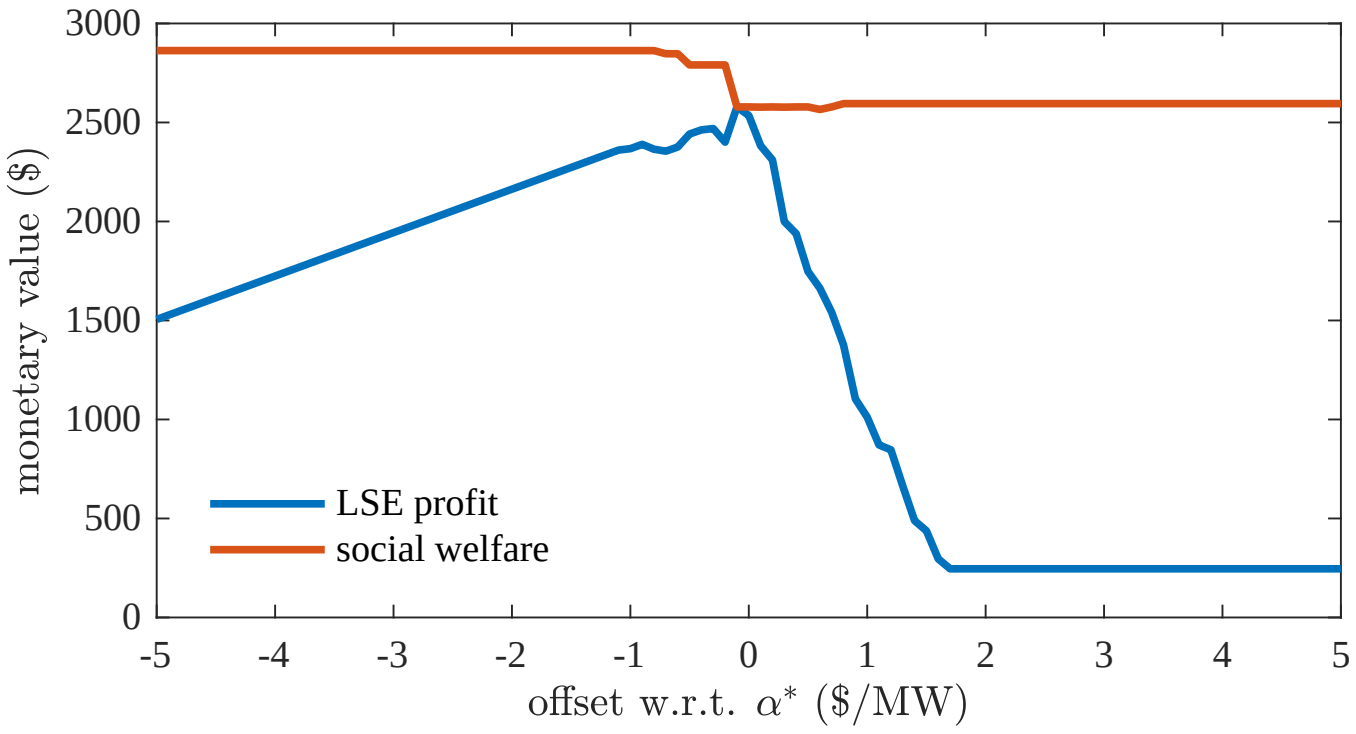

Fig. 2. WEM social welfare and LSE's profit versus LSE's energy price offset.

This is verified in Table I since the values 35.094, 35.325, and 36.299 exactly mathc the smallest value in $\boldsymbol{c}^{(k)}$.

We also contrast the base case with the case without joint bidding and pricing, where the LSE is non-profit and thus the market clearing price $\pi$ is posted directly to the EUCs. The results for this case are obtained by imposing an additional constraint in the MILP formulation that REM energy price $\alpha$ must equal to the LMP $\pi$ in the WEM, and the REM reserve price $\beta$ must equal to the reserve clearing price $\nu$ in the WEM. As shown in Table I, the profit-seeking LSE manages to manipulate the WEM to lower the clearing price for a greater difference from the optimal price it determines. In addition, a degradation of social welfare is expected due to this profit-seeking behavior of the LSE.

In addition to the non-profit case, other non-optimal pricing of the LSE may also have impacts on the clearing results. To verify this, we vary the value of $\boldsymbol{\alpha}$ by adding offsets with respect to (w.r.t.) the optimal one $\boldsymbol{\alpha}^*$. As shown in Fig. 2, the LSE's profit is indeed maximized at the optimal prices we proposed and changes almost monotonically before and after the peak point. Interestingly, the profit does not decrease to zero as $\alpha$ increases but rests at a constant level. We recognize that this level corresponds to the energy consumption level of $\overline{y}^{(k)}$, the upper bound of the reserve provided by the EUC. This is because the LSE seeks to prevent the energy demand from vanishing to zero (otherwise it will make no profit from energy), by encouraging the EUCs to provide as much reserve as possible, which serves as the lower bound for energy consumption (see (5a)). Moreover, the social welfare decreases with the increase of $\boldsymbol{\alpha}$ since a smaller demand of EUCs is satisfied.

### C. Impacts of Estimated Parameters

We also investigate how estimated parameters in the tri-layer model impact the clearing results. Towards this end, we first vary the parameters in the lower layer, e.g., the coefficients of EUCs' benefit functions. In particular, we vary $\boldsymbol{c}^{(k)}$ by adding offsets w.r.t. its values in the base case. As shown in Fig. 3, the LMP $\pi$ is always overlapped with the LSE's bids $a_j^d$ for $j = 1, 2, 3$, which correspond to the marginal cost of the power generation. This shows that one way for the profit-seeking LSE to manipulate the WEM is to submit bids that always equal to the marginal generation cost so that it can get its demand cleared at as low price as possible. This phenomena has also been observed in other simulation settings as shown later. Moreover, from Fig. 3, the increase of $\mathbf{c}^{(k)}$ in general causes the increase of LSE's optimal energy prices $\boldsymbol{\alpha}$. This is because the optimal price is essentially selected from the $\mathbf{c}^{(k)}$ values, i.e., the marginal benefit of some EUC's energy block.

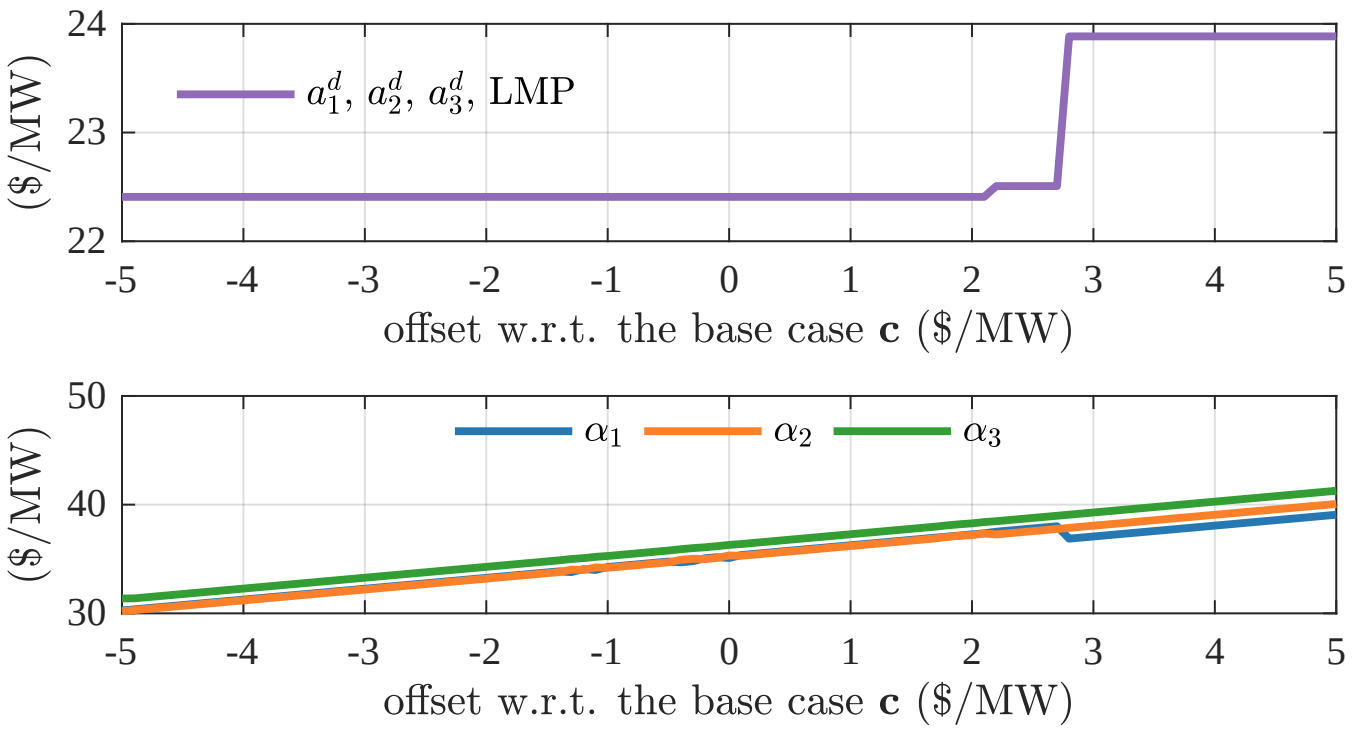

Fig. 3. LSE's bids $a_j^d$ and energy prices $\alpha_j^*$ for $j = 1, 2, 3$, LMP $\pi$ versus estimated EUCs' benefit coefficients $\boldsymbol{c}^{(k)}$ (with reserve).

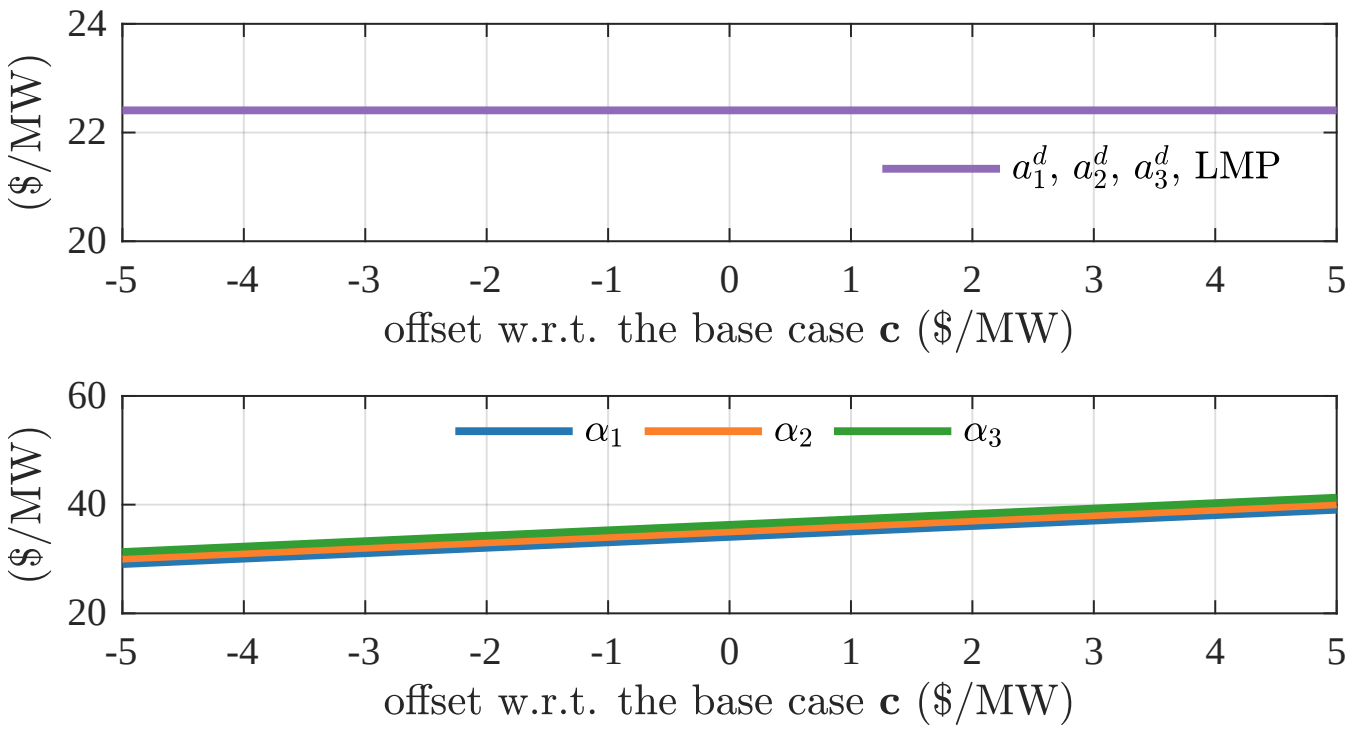

Fig. 4. LSE's bids $a_j^d$ and energy prices $\alpha_j^*$ for $j = 1, 2, 3$, LMP $\pi$ versus estimated EUCs' benefit coefficients $\boldsymbol{c}^{(k)}$ (without reserve).

Interestingly, however, this change is not monotone, which causes the abrupt switch of the marginal generator accordingly. We presume this is due to the presence of reserve trading in the market. In particular, as $\boldsymbol{\alpha}$ increases, the LMP $\pi$ is fixed as the marginal generation cost from Fig. 3. From the LSE's objective function, the first term of profit gained from the energy market will keep increasing and dominate, if the amount of energy demand $p_j^d$ is fixed. Thus, a little increase in the cleared amount $p_j^d$ can lead to a great increase in profit. Therefore, the LSE is willing to decrease the energy price a little to trade for the increase in amount, with the flexibility of making profit from the reserve market. With no reserve, however, any increase of $\boldsymbol{\alpha}$ that does not keep up with the increase of $\mathbf{c}^{(k)}$ (or even decrease) will result in a relatively low profit. To verify this, we simulate the setting without reserve and present the results in Fig. 4. It is seen that the energy price $\boldsymbol{\alpha}$ now increases monotonically with the marginal benefit $\mathbf{c}^{(k)}$ as expected, and the cleared energy demand does not change. This observation illustrates that the joint pricing and bidding problem becomes more involved when both reserve and energy are considered in our tri-layer market structure.

Furthermore, for the upper layer, we vary the estimated energy bid $a_j^d$ of other LSEs to evaluate how the possible inaccuracy of the estimation would impact the clearing results. As shown in Fig. 5, although other LSEs' bids increase, the LSE's bids still always align with the LMP $\pi$, ensuring that it can be cleared at a low LMP. We observe a peak in the LMP when the offset equals $-10$ ($/MW). This observation may be explained by the flexibility of the reserve in the market. In fact, as other LSEs' bids increases, the LMP tends to increase accordingly, which may decrease the LSE's profit from energy market. Therefore, the LSE may increase the energy price (as shown in Fig. 6) to reduce the amount of cleared energy so as to mitigate the profit loss resulted from it. Meanwhile, for compensation, the LSE can turn to the reserve market for

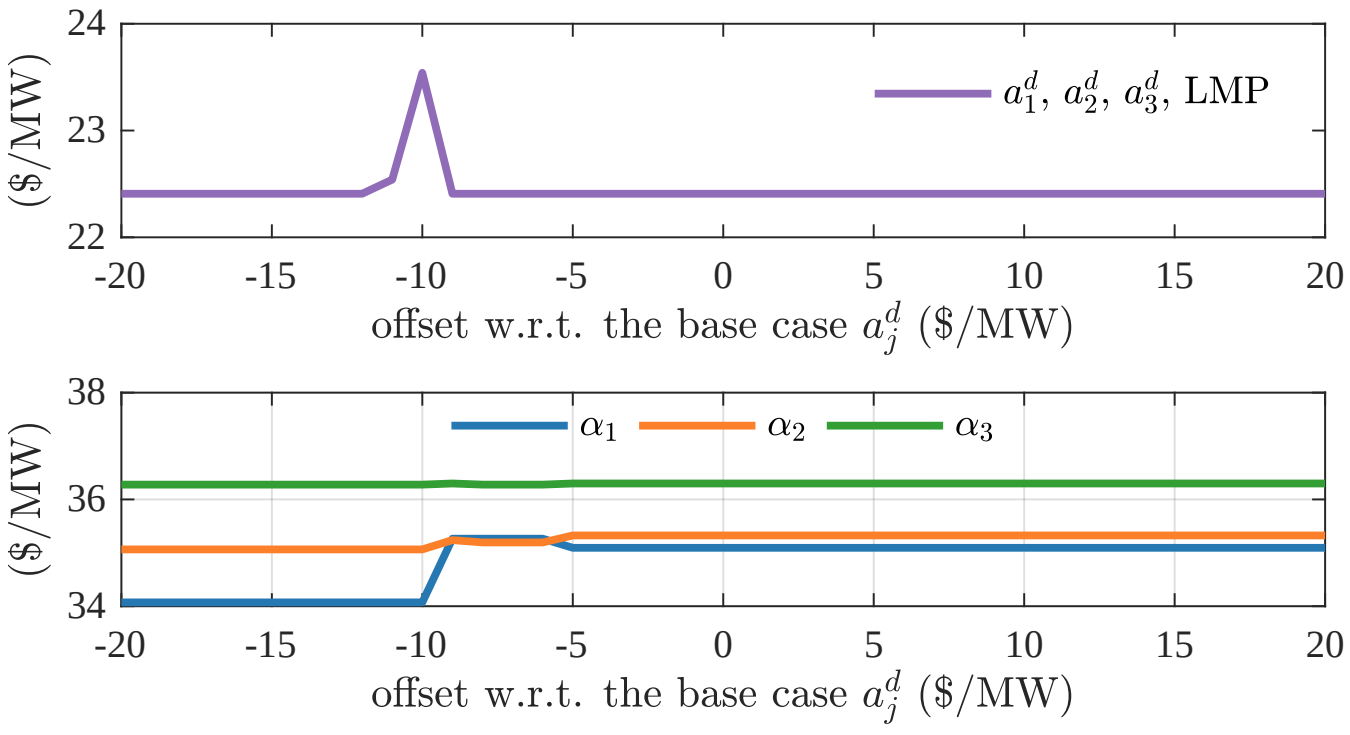

Fig. 5. LSE's bids $a_j^d$ and energy prices $\alpha_j^*$ for $j = 1, 2, 3$, LMP $\pi$ versus estimated offers of other LSEs $a_j^d$ for $j = 4, 5, 6$ (with reserve).

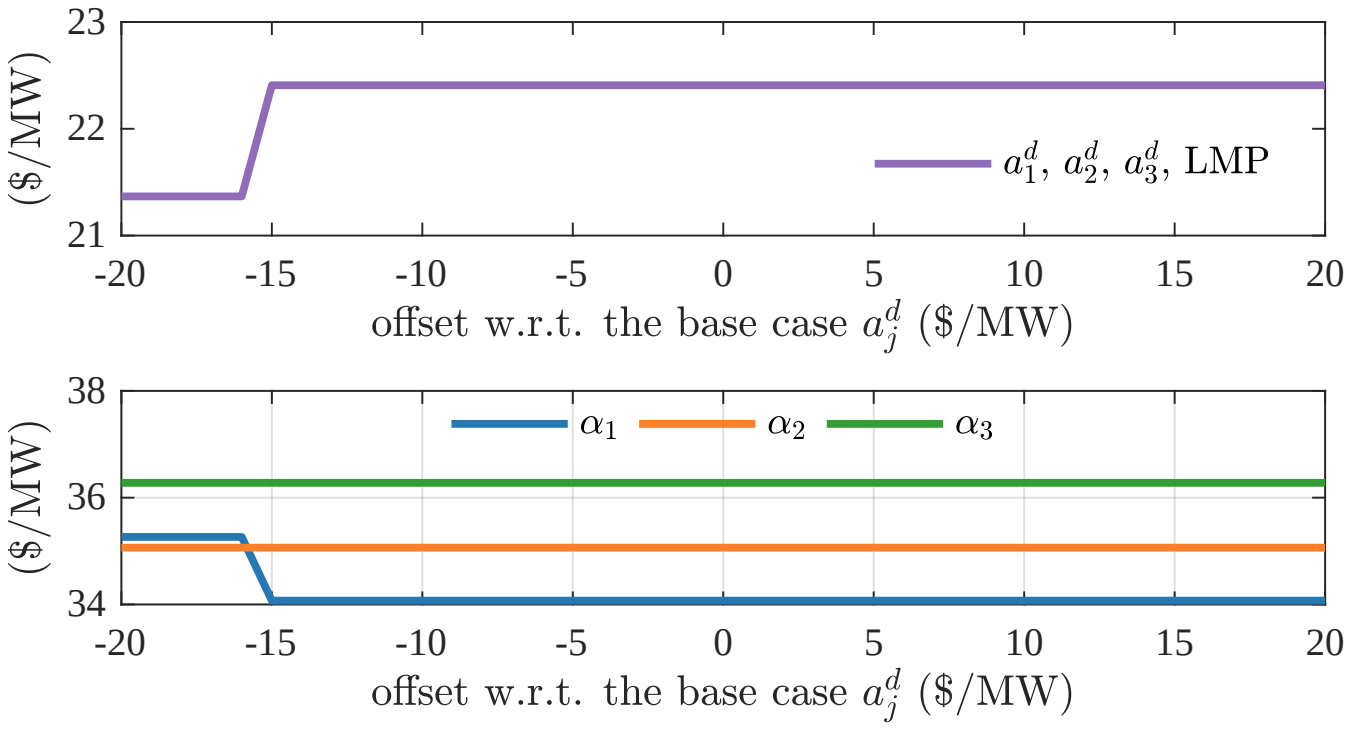

Fig. 6. LSE's bids $a_j^d$ and energy prices $\alpha_j^*$ for $j = 1, 2, 3$, LMP $\pi$ versus estimated offers of other LSEs $a_j^d$ for $j = 4, 5, 6$ (without reserve).

profit. This explanation is corroborated in Fig. 6 when no reserve is traded, where the LMP increases anyway and the LSE has no room to accommodate for profit but decreases the energy price to encourage energy consumption.

### D. *Impacts of Congestion*

To investigate the impact of line congestion on the LSE's joint bidding and pricing, we decrease the line active power capacity in the base case. As is shown in Table II, the bid $\mathbf{a}^d$ varies at different pricing areas (buses) $\mathcal{C}_k$ as well. Note that still, the clearing LMP $\pi^*$ is fully controlled by the bids $\mathbf{a}^d$. Interestingly, compared to the base case results in Table I, the LSE's profit decreases from \$2595.422 to \$965.958, together with the decrease of social welfare from \$2578.448 to \$1876.811. This is because in general fewer demands can be satisfied with limited transmission line active power capacity.

TABLE II
CLEARING RESULTS OF CASE WITH LINE CONGESTION AND CASE WITH EQUAL ENERGY PRICE AT DIFFERENT PRICING AREAS

| | $\alpha^*$ (\$/MW) | | | profit (\$) |
|---|---|---|---|---|
| | $\mathcal{C}_1$ | $\mathcal{C}_2$ | $\mathcal{C}_3$ | |
| congested | 35.914 | 36.919 | 37.858 | 965.958 |
| equal price | 37.834 | 37.834 | 37.834 | 753.542 |
| | $\boldsymbol{a}^{d*}$ (\$/MW) | | | welfare (\$) |
| | $\mathcal{C}_1$ | $\mathcal{C}_2$ | $\mathcal{C}_3$ | |
| congested | 23.885 | 23.275 | 22.408 | 1876.811 |
| equal price | 23.885 | 23.223 | 22.750 | 1901.213 |
| | LMP $\pi^*$ (\$/MW) | | | |
| | $\mathcal{C}_1$ | $\mathcal{C}_2$ | $\mathcal{C}_3$ | |
| congested | 23.885 | 23.275 | 22.408 | |
| equal price | 23.885 | 23.223 | 22.750 | |

In practical electricity market, one LSE may not be allowed to design different energy and reserve prices at its pricing areas. Consider this regulation constraint, the LSE may not be able to achieve the optimal profit. This is corroborated in Table II, where second rows in the table correspond to the case with $\alpha_1 = \alpha_2 = \alpha_3$. The resulted profit \$753.542 is less than the optimal \$965.958.

### E. Large System Test Case

We further test the proposed model using the IEEE 118-bus test case [26]. We randomly choose 6 out of the total 64 load buses, buses $2, 11, 16, 35, 47, 83$, to represent different pricing areas, each aggregating 10 EUCs, with each EUC submitting 3 bids as before. The benefit coefficients $\boldsymbol{c}$ and $\boldsymbol{d}$ are i.i.d. samples drawn uniformly from $[34, 38]$ and $[5, 7]$, respectively. The marginal generation cost $a_j^g$ and $b_j^g$ are sampled uniformly from $[21, 25]$ and $[4.5, 7.5]$, respectively. Other LSEs' bids $a_j^d$ and offers $b_j^d$ are sampled uniformly from $[33, 39]$ and $[5.5, 8.5]$, respectively. The upper limits $\bar{\boldsymbol{x}}^{(k)}$ and $\bar{\boldsymbol{y}}^{(k)}$ are selected as in the small system case. The line active power limits are selected randomly from the range of $[100, 300]$ MW.

Similar to the clearing results in the base case, the LSE is able to pursue profits by designing various energy and reserve prices $\boldsymbol{\alpha}^*$ and $\boldsymbol{\beta}^*$ across pricing areas, instead of posting the clearing price $\pi$ to the EUCs directly, as is shown in Table III. Note that since line congestions are considered, the resulted LMP for both with and without LSE's joint pricing and bidding vary across different pricing areas. Moreover, in the large system case, the LSE's profit-seeking behavior does not vary the social welfare significantly from the case without joint biding/pricing. This suggests that the LSE has smaller market power in this larger system. In fact, it indeed controls fewer pricing areas (buses) (6 out of 64) than in the small system case (3 out of 6).

We also present the computation time and MILP gaps for solving both cases in Table IV. It is seen that the proposed approach can efficiently solve the joint pricing and bidding problem with high precision. Moreover, the proposed methodology is scalable to the large system test case.

TABLE III
LARGE SYSTEM CLEARING RESULTS WITH/W.O. JOINT BIDDING/PRICING

| | $\boldsymbol{\alpha}^*$ (\$/MW) | | | | | |
|---|---|---|---|---|---|---|
| | $\mathcal{C}_1$ | $\mathcal{C}_2$ | $\mathcal{C}_3$ | $\mathcal{C}_4$ | $\mathcal{C}_5$ | $\mathcal{C}_6$ |
| with | 34.109 | 34.073 | 37.003 | 35.307 | 34.050 | 37.954 |
| w.o. | 21.928 | 31.954 | 35.400 | 34.689 | 21.711 | 37.880 |
| | $\boldsymbol{a}^{d*}$ (\$/MW) | | | | | |
| | $\mathcal{C}_1$ | $\mathcal{C}_2$ | $\mathcal{C}_3$ | $\mathcal{C}_4$ | $\mathcal{C}_5$ | $\mathcal{C}_6$ |
| with | 21.928 | 30.791 | 22.435 | 24.646 | 21.434 | 40.873 |
| w.o. | 21.928 | 31.954 | 35.400 | 34.689 | 21.711 | 37.880 |
| | $\boldsymbol{\beta}^*$ (\$/MW) | | | | | |
| | $\mathcal{C}_1$ | $\mathcal{C}_2$ | $\mathcal{C}_3$ | $\mathcal{C}_4$ | $\mathcal{C}_5$ | $\mathcal{C}_6$ |
| with | 5.527 | 5.368 | 5.645 | 6.655 | 5.514 | 5.872 |
| w.o. | 6.061 | | | | | |
| | LMP $\pi^*$ (\$/MW) | | | | | |
| | $\mathcal{C}_1$ | $\mathcal{C}_2$ | $\mathcal{C}_3$ | $\mathcal{C}_4$ | $\mathcal{C}_5$ | $\mathcal{C}_6$ |
| with | 21.928 | 30.791 | 22.435 | 24.646 | 21.434 | 40.873 |
| w.o. | 21.928 | 31.954 | 35.400 | 34.689 | 21.711 | 37.880 |
| | profit (\$) | | | welfare (\$) | | |
| with | 1116.722 | | | 18837.373 | | |
| w.o. | 0 | | | 19191.994 | | |

TABLE IV
COMPUTATION TIME AND MILP GAPS

| | presolve time (s) | computation time (s) | MILP gap (%) |
|---|---|---|---|
| 9-bus | 0.02 | 2.17 | 0.45 |
| 118-bus | 0.41 | 36.34 | 0.29 |

## VI. CONCLUDING REMARKS

In this paper, we proposed a formulation for the optimal joint bidding and pricing problem of a profit-seeking LSE based on a tri-layer market model. We show the joint bidding and pricing problem can be formulated as an MILP problem by leveraging the duality theory. Consequently, despite the complex structure of the model, the problem can still be solved efficiently, making it applicable for practical applications. We corroborated the applicability of the proposed model through numerical simulations, and revealed that the LSE's optimal strategy is influenced by the EUCs' benefit functions, other LSEs' biding/pricing strategies, and network congestions.

Future works will focus on developing a model with general-form energy benefit functions and reserve cost function of the EUCs. Besides, estimation methods for unknown parameters such as the coefficients of energy benefit function and reserve cost function of the EUCs, as well as the bids/offers from GENs and other LSEs will also be investigated.

**Hanchen Xu** (S'12-) received his B.Eng. and M.S. in Electrical Engineering from Tsinghua University in 2012 and 2014, respectively. He also obtained a M.S. in Applied Mathematics from University of Illinois at Urbana-Champaign in 2017. He is now a Ph.D. candidate in the Department of Electrical and Computer Engineering, University of Illinois at Urbana-Champaign. His current research interests include control, optimization, machine learning, with applications to smart grids and electricity markets.

**Kaiqing Zhang** (S'16-) received B.Eng. degree from Tsinghua University, Beijing, China in 2015. He is pursuing the Ph.D. degree in at the Department of Electrical and Computer Engineering, University of Illinois at Urbana-Champaign. His current research interests include reinforcement learning, optimization, and game theory, and their applications toward multi-agent systems.

**Junbo Zhang** (S'10-M'14) received his B.Eng. and Ph.D. from Tsinghua University in 2008 and 2013, respectively, followed by a postdoc fellow at the same university. He studied at The Hong Kong Polytechnic University from 2009 to 2011 under a collaborative program between Tsinghua University and The Hong Kong Polytechnic University. He is now an Associate Professor and Dean Assistant at the School of Electric Power, South China University of Technology. His research interests include stability, control, and knowledge automation of power systems, and electricity markets.